\input amstex
\magnification=1200
\documentstyle{amsppt}
\NoRunningHeads
\NoBlackBoxes
\define\Diff{\operatorname{Diff}}
\define\BU{\bold U}
\define\BB{\bold B}
\define\BE{\bold E}
\define\BP{\bold P}
\define\CM{\Cal M}
\define\CO{\Cal O}
\define\Fg{\frak g}
\define\HS{\Cal H\Cal S}
\define\sltwo{\operatorname{\frak s\frak l}(2,\Bbb C)}
\define\Vect{\operatorname{Vect}}
\define\CVect{\operatorname{\Bbb CVect}}
\define\Ner{\operatorname{Ner}}
\define\Rot{\operatorname{Rot}}
\define\form{\operatorname{form}}
\define\spn{\operatorname{span}}
\define\End{\operatorname{End}}
\topmatter
\title On the infinite--dimensional hidden symmetries. I.
Infinite dimensional geometry of $q_R$--conformal symmetries.
\endtitle
\author Denis V. Juriev
\endauthor
\affil\rm
``Thalassa Aitheria'' Research Center for Mathematical Physics
and Informatics,\linebreak
ul.Miklukho-Maklaya 20-180, Moscow 117437 Russia.\linebreak
E-mail: denis\@juriev.msk.ru\linebreak
\ \linebreak
November, 28, 1996\linebreak
E-version: funct-an/9612004
\endaffil
\abstract The infinite dimensional geometry of the $q_R$--conformal
symmetries (the hidden symmetries in the Verma modules over $\sltwo$
generated by the spin 2 tensor operators) is discussed.
\endabstract
\endtopmatter
\document
This paper opens the series of articles supplemental to [1], which also lie
in lines of the general ideology exposed in the review [2]. The main purpose
of further activity, which has its origin and motivation presumably in the
author's applied researches [3] on the interactively controlled systems
(i.e. the controlled systems, in which the control is coupled with unknown
or uncompletely known feedbacks), is to explicate the essentially {\sl
infinite--dimensional\/} aspects of the hidden symmetries, which appear in the
representation theory of the finite dimensional Lie algebras and related
algebraic structures. The series is organized as a sequence of topics, which
illustarate this basic idea on the simple and tame examples without superfluous
difficulties and details as well as in the series [1]. Many of objects, which
will appear, are somehow related to ones discussed previously. However, the
matherial will be treated more geometrically, presumably, from the points of
view of the infinite dimensional geometry (cf.[4]), an infinite dimensional
version of the nonlinear geometric algebra (cf.[5]) and the infinite
dimensional noncommutative geometry (cf.[6]). The intent reader will see a lot
of ideological similarities between the subject of the paper and the infinite
dimensional geometric picture for a second quantized string (see [7] and
numerous refs wherein).

\head 1. Unitary $\HS$--pseudorepresentations of $\widetilde\Diff_+(S^1)$,
infinite dimensional geometry of $q_R$--conformal symmetries and the related
topics\endhead

\definition{Definition 1A} {\it A unitary $\HS$--pseudorepresentation\/} of the
Lie group $G$ in the Hilbert space $H$ is a homomorphism of $G$ into
$\BU(H)/\BU_0(H)$, where $\BU(H)$ is the group of all unitary operators in $H$
and $\BU_0(H)$ is its subgroup of operators $U$ such that $1-U$ is a
Hilbert--Schmidt operator.
\enddefinition

The motivation for the choice of the name ``$\HS$--pseudorepresentation'' is an
analogy with the constructions of the pseudodifferential calculus [8], the
abbreviation $\HS$ means that an admissible deviation from the ``true''
representations belongs to the class $\HS$ of the Hilbert--Schmidt operators
(one may also consider the operator ideals of compact or trace-class operators
instead of $\HS$, it allows to generalize the constructions to the linear
topological spaces). It is rather convenient to consider the pullbacks
$\tilde\pi:G\mapsto\BU(H)$ instead of the $\HS$--pseudorepresentations $\pi$
themselves. The pullback is called {\it monoassociative\/} iff
$\tilde\pi(g_{t+s})=\tilde\pi(g_t)\tilde\pi(g_s)$ for any one--parameter
subgroup $\{g_t\}$ of $G$. The pillback $\tilde\pi$ obeys the following claim:
$$\tilde\pi(g_1)\tilde\pi(g_2)\equiv\tilde\pi(g_1g_2)\pmod{\HS}$$
for $g_1$, $g_2$ from $G$.

\definition{Definition 1B} Two unitary $\HS$--pseudorepresentations $\pi_1$
and $\pi_2$ of the Lie group $G$ in the Hilbert spaces $H_1$ and $H_2$ are
called {\it equivalent\/} iff there exists an isomorphism $S:H_1\mapsto H_2$
of the linear spaces such that for any element $g$ of the group $G$ there
exists an element $U_g$ of $\BU_0(H)$ such that
$$\tilde\pi_2(g)=U_gS\tilde\pi_1(g)S^{-1},$$
where $\tilde\pi_i$ are the mappings of $G$ into $\BU(H_i)$, which are the
pullbacks of the homomorphisms $\pi_i$ alongside the natural projections
$\varepsilon_i:\BU(H_i)\mapsto\BU(H_i)/\BU_0(H_i)$.
\enddefinition

The infinitesimal counterpart of the notion of the unitary
$\HS$--pseu\-do\-rep\-re\-sen\-tation of the Lie group $G$ is one of the
$\HS$--projective representation of the Lie algebra $\Fg$ [1:Topic 10]. The
idea of the unitary $\HS$--pseudorepresentation is ideologically very closed
to one of the quasirepresentation of A.I.Shtern [9] (A.I.Shtern also defined
the pseudorepresentations in a way analogous to the definition of a
pseudocharacter, however, his definition differs from ours. Below we shall
not discuss A.I.Shtern's pseudorepresentations).

Let us now consider the $\HS$--projective representations $\tau_h$ of the Witt
algebra in the unitarizable Verma modules $V_h$ ($h$ is the highest weight)
over the Lie algebra $\sltwo$. Such representations realize the generators
of the Witt algebra as tensor operators of spin 2 in the Verma modules or,
otherwords, as the $q_R$--conformal symmetries [3]. The $\HS$--projective
representations $\tau_h$ may be extended to the $\HS$--projective
representations $\tau_h$ of the algebra $\Vect(S^1)$ of the smooth vector
fields on a circle in the suitable G\aa rding space. Let
$\widetilde\Diff_+(S^1)$ be the universal covering of the group $\Diff_+(S^1)$
of the orientation preserving diffeomorphisms of a circle $S^1$.

\proclaim{Theorem 1} The $\HS$--projective representations $\tau_h$ of
the Lie algebra $\Vect(S^1)$ are exponentiated to the unitary
$\HS$--pseudorepresentations $T_h$ of $\widetilde\Diff_+(S^1)$ in
the Hilbert spaces $H_h$, which are the completions of the unitarizable
Verma modules $V_h$ over the Lie algebra $\sltwo$.
\endproclaim

\remark{Remark 1} There exists a uniquely defined continuous monoassociative
pullback $\tilde T_h$ of the unitary $\HS$--pseudorepresentation $T_h$ such
that the infinitesimal operators for the generators of $\Vect(S^1)$ coincide
with one of the $\HS$--projective representation of the Witt algebra in the
Verma module over the Lie algebra $\sltwo$ by the tensor operators of spin 2
[1:Topic 10].
\endremark

\remark{Remark 2} The theorem may be generalized to the non-unitarizable
Verma modules if one changes the class $\HS$ of the Hilbert-Schmidt operators
to the class of all compact operators.
\endremark

Let $\CM=\Diff_+(S^1)/S^1$ be the flag manifold for the Virasoro--Bott
group [10;4].

\remark{Corollary} The monoassociative pullback $\tilde T_h$ defines an
imbedding $\epsilon:\CM\mapsto\BP(H_h)$ of the infinite--dimensional
manifold $\CM$ into the projective space $\BP(H_h)$ as an orbit $\CO$ of the
highest vector.
\endremark

Let us consider the real projective space $\BP(H_h^{\Bbb R})$ as a fiber
bundle over $\BP(H_h)$ with fibers being isomorphic to $S^1$ and the projection
$p:\BP(H_h^{\Bbb R})\mapsto\BP(H_h)$. The subset $\tilde\CO=p^{-1}(\CO)$ of
$\BP(H_h^{\Bbb R})$ may be identified with $\Diff_+(S^1)$. Let us fix any point
$x$ of $\tilde\CO$. The operators $\tilde T_h(g)$ ($g\in\Diff_+(S^1)$) in $H_h$
may be considered as operators in $\BP(H_h^{\Bbb R})$. The orbit of $x$ under
these operators supply $\tilde\CO$ by another structure of the groupuscular
(local group) centered in the point $x$. The obtained groupuscular structure
on $\tilde\CO\simeq\Diff_+(S_1)$ is not canonical (see [5]), however, the
operators of deviation of the new structure from the canonical belong to
the Hilbert-Schmidt class $\HS$. Note though the pullback $\tilde T_h$ is
monoassociative the noncanonical groupuscular structure on $\Diff_+(S_1)$
is not odular (see [5]) as an artefact of the infinite dimensionality (the
image of the exponential map of $\Vect(S^1)$ does not cover any neighbourhood
of the identity in $\Diff_+(S^1)$).

\remark{Remark 3} The noncanonical groupuscular structure on $\Diff_+(S^1)$
realizes this group as a transformation quasigroup [11] on the orbit $\CO$ (a
transformation pseudogroup in terms of [12]).
\endremark

Let us now discuss the complex structure on the orbit $\CO$.

\definition{Definition 2} {\it A linear $\HS$--pseudorepresentation\/} of the
Lie group $G$ in the Hilbert space $H$ is a homomorphism of $G$ into
$\BB(H)/\BB_0(H)$, where $\BB(H)$ is the group of all invertible bounded
operators in $H$ and $\BB_0(H)$ is its subgroup of operators $A$ such that
$1-A$ is a Hilbert--Schmidt operator. Two linear representations $\pi_1$ and
$\pi_2$ of the Lie group $G$ in the Hilbert spaces $H_1$ and $H_2$ are called
{\it equivalent\/} iff there exists an isomorphism $S:H_1\mapsto H_2$ of
the linear spaces such that for any element $g$ of the group $G$ there exists
an element $A_g$ of $\BB_0(H)$ such that
$$\tilde\pi_2(g)=A_gS\tilde\pi_1(g)S^{-1},$$
where $\tilde\pi_i$ are the mappings of $G$ into $\BB(H_i)$, which are the
pullbacks of the homomorphisms $\pi_i$ alongside the natural projections
$\varepsilon_i:\BB(H_i)\mapsto\BB(H_i)/\BB_0(H_i)$.
\enddefinition

This definition may be generalized on any topological semigroup $\Gamma$
if one change the group $\BB(H)$ of all invertible bounded operators to the
semigroup $\BE(H)$ of all bounded operators in $H$.

Let $\widetilde\Ner$ be the universal covering of the Neretin semigroup, the
mantle of the group $\Diff_+(S^1)$ [13]. The semigroup $\widetilde\Ner$ is a
complex semigroup so one may consider its holomorphic linear
$\HS$--pseudorepresentations.

\proclaim{Theorem 2} The unitary $\HS$--pseudorepresentations $T_h$ of
$\widetilde\Diff_+(S^1)$ in the Hil\-bert spaces $H_h$ may be extended to the
holomorphic linear $\HS$--pseu\-do\-rep\-re\-sen\-tations of the Neretin
semigroup $\Ner$.
\endproclaim

\remark{Remark 4} The theorem 4 may be generalized to the non-unitarizable
Verma modules if one changes the class $\HS$ of the Hilbert-Schmidt operators
to the class of all compact operators.
\endremark

\remark{Remark 5} The imbedding $\varepsilon:\CM\mapsto\BP(H_h)$ is
holomorphic.
\endremark

Remark 5 supplies us by a new construction of the Kirillov's complex structure
on $\CM=\Diff_+(S^1)/S^1$.

Note that the infinite--dimensional complex manifold $\CM$ may be identified
with the class $S_0$ of univalent functions [10] (see also [4,14]). The
imbedding $\varepsilon:S_0\mapsto\BP(H_h)$ may be analytically extended to
the imbedding $\tilde\varepsilon$ of the space $\Bbb C_0[[z]]$ of all formal
power series of the form $z+c_1z^2+c_2z^3+\ldots+c_nz^{n+1}+\ldots$ into the
projective space $\BP(V_h^{\form})$ over the formal Verma module $V_h^{\form}$
for the Lie algebra $\sltwo$ (cf.[14]). Any formal Verma module over $\sltwo$
is identified with the space of the formal power series $\Bbb C[[z]]$ whereas
the Verma module itself is identified with the space of polynomials $\Bbb
C[z]$.

\remark{Remark 6 {\rm (The hypothetical criterion of the univalence
(cf.[14]))}}
$$\tilde\epsilon^{-1}(\tilde\epsilon(\Bbb C_0[[z]])\cap\BP(H_h))=S_0$$
\endremark

\remark{Corollary {\rm (hypothetical)}} Let $\deg(c_n)=n$, then there exists
a universal sequence $\{P_m\}$ of the polynomials $P_m=P_m(c_1,\ldots c_m)$,
$\deg(P_m)=m$ such that the formal power series
$z+c_1z^2+c_2z^3+\ldots+c_nz^{n+1}+\ldots$ defines an univalent function
$f(z)$ from the class $S_0$ iff
$\sum_{m=1}^{\infty}|P_m(c_1,\ldots c_m)|^2<\infty$.
\endremark

Note that it follows from the results of [14] that there exists a universal
sequence $\{P_m^{(k)}\}$ ($1\le k\le p(m)$, $p(m)$ is the partition function)
of polynomials $P_m^{(k)}=P_m^{(k)}(c_1,\ldots c_m)$, $\deg(P_m^{(k)})=m$ such
that the formal power series $z+c_1z^2+c_2z^3+\ldots+c_nz^{n+1}+\ldots$ defines
an univalent function $f(z)$ from the class $S_0$ iff
$\sum_{m=1}^{\infty}\sum_{k=1}^{p(m)}|P_m^{(k)}(c_1,\ldots c_m)|^2<\infty$.
The statement of the corollary is essentially stronger.

\head 2. Composed representations of the Witt isotopic pair from
$q_R$--conformal stress-energy tensor and the related ($q_R$--affine)
current
\endhead

\definition{Definition 3A {\rm [15] (see also [2:\S2.2;16])}} The pair
$(V_1,V_2)$ of linear spaces is called {\it an (even) isotopic pair\/} iff
there are defined two mappings $m_1:V_2\otimes\bigwedge^2V_1\mapsto V_1$ and
$m_2:V_1\otimes\bigwedge^2V_2\mapsto V_2$ such that the mappings $(X,Y)\mapsto
[X,Y]_A=m_1(A,X,Y)$ ($X,Y\in V_1$, $A\in V_2$) and $(A,B)\mapsto
[A,B]_X=m_2(X,A,B)$ ($A,B\in V_2$, $Y\in V_1$) obey the Jacobi identity for
all values of a subscript parameter (such operations will be called {\it
isocommutators\/} and the subscript parameters will be called {\it isotopic
elements\/}) and are compatible to each other, i.e. the identities
$$\align
[X,Y]_{[A,B]_Z}=&\tfrac12([[X,Z]_A,Y]_B+[[X,Y]_A,Z]_B+[[Z,Y]_A,X]_B-\\
&[[X,Z]_B,Y]_A-[[X,Y]_B,Z]_A-[[Z,Y]_B,X]_A)\endalign$$
and
$$\align
[A,B]_{[X,Y]_C}=&\tfrac12([[A,C]_X,B]_Y+[[A,B]_X,C]_Y+[[C,B]_X,A]_Y-\\
&[[A,C]_Y,B]_X-[[A,B]_Y,C]_X-[[C,B]_Y,A]_X)\endalign$$
($X,Y,Z\in V_1$,
$A,B,C\in V_2$) hold.
\enddefinition

The even isotopic pairs are just the anti-Jordan pairs of J.R.Faulkner and
J.C.Fer\-rar [17] if the characteristics of the basic field is not equal to
2 (see [15;2:\S2.2]).

\remark{Remark 7} Let $M$ be an arbitrary (smooth) manifold, then the space
$\CO(M)$ of all smooth functions on $M$ and the space $\Vect(M)$ of all
smooth vector fields on $M$ form an isotopic pair with isocommutators
$$[v_1,v_2]_f=\Cal L_{v_1}(f)v_2-\Cal L_{v_2}(f)v_1+f[v_1,v_2]\quad
(v_1,v_2\in\Vect(M),\ f\in\CO(M))$$
and
$$[f_1,f_2]_v=\Cal L_v(f_2)f_1-\Cal L_v(f_1)f_2\quad (f_1,f_2\in\CO(M),
\ v\in\Vect(M)),$$
here $\Cal L$ denotes the Lie derivative. The constructed isotopic pair is
called {\it the geometric isotopic pair\/} [15].
\endremark

Let us consider the complexification of the geometric isotopic pair with
$M=S^1$. Let $e_k=ie^{ikt}\partial_t$ and $f_k=ie^{ikt}$ be the natural
basises in the spaces $\CO^{\Bbb C}(S^1)$ and $\CVect(S^1)$, $t$ is an angle
coordinate on the circle $S^1$. The isocommutators have a nice symmetric
form [15]:
$$[e_i,e_j]_{f_k}=(i-j)e_{i+j+k},\qquad [f_i,f_j]_{e_k}=(i-j)f_{i+j+k}.$$
The isotopic pair $(V_1,V_2)$ formally generated by the elements $e_k$
($k\in\Bbb Z$) and $f_k$ ($k\in\Bbb Z$) will be called {\it the Witt
isotopic pair}.

The isocommutators in $(V_1,V_2)$ are r-matrix, i.e. they may be received
from the canonical brackets in the Witt algebra by use of the classical
r-matrices [18]. The brackets in the Witt algebra have the form
$[e_i,e_j]=(i-j)e_{i+j}$ (and similar for $f_k$), the classical r-matrices
$R_x$ ($[a,b]_x\!=\![R_x(a),b]\!+\![a,R_x(b)]$) have the form
$R_{f_k}(e_i)=e_{i+k}$ and $R_{e_k}(f_i)=f_{i+k}$. These r-matrices do not
obey the modified Yang--Baxter equation (see [18]). Note that the mapping
$x\mapsto R_x$ is multiplicative, i.e. $R_{f_i}R_{f_j}=R_{f_{i+j}}$ and
$R_{e_i}R_{e_j}=R_{e_{i+j}}$.

\definition{Definition 3B {\rm [15]}} The pair $(V_1,V_2)$ of linear spaces is called
{\it the (even) isotopic composite\/} iff there is fixed a set of subpairs
$(V_1^\alpha,V_2^\alpha)$ of the pair $(V_1,V_2)$ such that each pair
$(V_1^\alpha,V_2^\alpha)$ is supplied by the structure of an (even) isotopic
pair. The structures of the isotopic pairs are compatible, it means that the
restrictions of the structures for two different $\alpha$ and $\beta$ on the
pair $(V_1^\alpha\cap V_1^\beta,V_2^\alpha\cap V_2^\beta)$ coincide.
The isotopic composite is called {\it dense\/} iff $\biguplus_\alpha
V_1^\alpha=V_1$ and $\biguplus_\alpha V_2^\alpha=V_2$ (here $\biguplus$ denotes
the sum of subspaces). The Lie composite is called {\it connected\/}
iff for all $\alpha$ and $\beta$ there exists a sequence $\gamma_1,\ldots
\gamma_m$ ($\gamma_1=\alpha$, $\gamma_m=\beta$) such that
$(V_1^{\gamma_l},V_2^{\gamma_l})\cap(V_1^{\gamma_{l+1}},V_2^{\gamma_{l+1}})
\ne\varnothing$.
\enddefinition

\example{Example {\rm (The Tetrahedron Isotopic Composite)}}
Let us consider a tetrahedron with vertices $A$, $B$, $C$, $D$ and the
faces $X=(BCD)$, $Y=(CDA)$, $U=(DBA)$, $V=(BCA)$. The linear spaces
$V_1$ and $V_2$ are spanned by the generators $e_A$, $e_B$, $e_C$, $e_D$,
$f_X$, $f_Y$, $f_U$, $f_V$ labelled by the vertices and the faces.
The subpairs $(V_1^\alpha,V_2^\alpha)$ are spanned by the generators
corresponding either to three vertices and a face between them or
by three faces and their common vertex. The isotopic pairs
$(V_1^\alpha,V_2^\alpha)$ are isotopic pairs associated with the Lie
algebra $\operatorname{so}(3)$ [15,16]. The isocommutators are compatible
with the orientation.
\endexample

The Witt isotopic pair $(V_1,V_2)$ is supplied by the canonical structure
of the dense and connected isotopic composite: $V_1^1=\spn(e_i, i\ge-1)$,
$V_2^1=\spn(f_i, i\ge 0)$, $V_1^2=\spn(e_i, i\le 1)$, $V_2^2=\spn(f_i, i\le
0)$.

\definition{Definition 4 {\rm [15] (see also [2:\S2.2])}}

{\rm A.} {\it A representation\/} of the (even) isotopic pair $(V_1,V_2)$ in
the linear space $H$ is the pair of linear mappings $(T_1,T_2)$ from the
spaces $V_1$ and $V_2$, respectively, to $\End(H)$ such that
$$\aligned
T_1([X,Y]_A)&=T_1(X)T_2(A)T_1(Y)-T_1(Y)T_2(A)T_1(X),\\
T_2([A,B]_X)&=T_2(A)T_1(X)T_2(B)-T_2(B)T_1(X)T_2(A),
\endaligned
$$
where $X,Y\in V_1$ and $A,B\in V_2$.

{\rm B.} {\it A representation\/} of the (even) isotopic composite $(V_1,V_2)$
in the linear space $H$ is the set of representations
$(T_1^\alpha,T_2^\alpha)$ of the isotopic pairs $(V_1^\alpha,V_2^\alpha)$
compatible to each other. The compatibility means that the restrictions of
the representations for different $\alpha$ and $\beta$ on the intersections
of the corresponding isotopic pairs coincide.

{\rm C.} Let an isotopic pair $(V_1,V_2)$ is supplied by the structure of
the isotopic composite (by the fixing of a set of its proper subpairs), then
the representations of the least will be called {\it the composed
representations\/} of $(V_1,V_2)$.
\enddefinition

\remark{Remark 8} All the objects from the definition 4 have their higher
combinatorial generalizations ({\it graph--representations\/}) -- cf.[15].
\endremark

\remark{Remark 9} The representations of the tetrahedron isotopic composite
realize the composed representations of the isoquaternionic isotopic pair [15].
\endremark
Note that the Witt algebra admits a composed representation [1:Topic 9]
in the Verma module $V_h$ over the Lie algebra $\sltwo$ by the tensor
operators of spin 2 (i.e. the $q_R$--conformal symmetries).

\proclaim{Theorem 3} The composed representation of the Witt algebra in the
Verma module $V_h$ over the Lie algebra $\sltwo$ may be extended to the
composed representation of the Witt isotopic pair by the tensor operators
of spin 1 and 2, namely
$$\aligned
e_k\mapsto(\xi+(k+1)h)\partial_z^k\quad (k\ge 0),&\qquad
e_{-k}\mapsto z^k\tfrac{\xi+(k+1)h}{(\xi+2h)\ldots(\xi+2h+k-1)}\quad
(k\ge 1),\\
f_k\mapsto\partial_z^k\quad (k\ge 0),&\qquad
f_{-k}\mapsto z^k\tfrac1{(\xi+2h)\ldots(\xi+2h+k-1)}\quad (k\ge 1),
\endaligned
$$
where $\xi=z\partial_z$.
\endproclaim

Note that operators $T_2^1(f_1)$ and $T_2^2(f_{-1})$ generate the
Lobachevski{\v\i}--Berezin algebra (see e.g.[1,3]) whereas $T_1^1(e_k)$,
$T_1^2(e_{-k})$ ($k\ge 2$) and $T_1^1(e_0)=T_1^2(e_0)$ generate the nonlinear
$\operatorname{sl}_2$ in sense of M.Ro\v cek [19].

The generating function
for the tensor operators of spin 2 is the $q_R$-conformal stress-energy
tensor whereas one for the tensor operators of spin 1 is the $q_R$-affine
current. The $q_R$-conformal stress-energy tensor may be received from
the $q_R$-affine current (which is involutive [3]) as by the truncated
Sugawara construction as an exponential of the associated Fubini-Veneziano
field [3].

\remark\nofrills{Problems:}
\roster
\item"--" To globalize the Witt isotopic pair, i.e. to construct the
object, whose infinitesimal counterpart is just the Witt isotopic pair.
There are hopes that the natural framework for the globalization may be found
in lines of [20] (``the universal quantum group''). The result should be
described as a quantum group with the multiple noncommutative parameters
$q_i$\footnote"*"{\ The idea to consider the multiparametric quantum
deformations was explored by Yu.I.Manin and his pupils as I know. However,
the multiple parameters of their deformations always commute.\newline}, which
form also a quantum group with elements of the first one as the multiple
noncommutative parameters of quantization.
\item"--" To globalize the construction of the composed representation of
the Witt composite in the Verma modules $V_h$ over the Lie algebra $\sltwo$.
\item"--" To describe an algebraic structure analogous to the Witt
isotopic composite and represented by tensor operators of the higher spins.
\endroster
\endremark

\Refs
\roster
\item" [1]" Juriev D., Topics in hidden symmetries. I--V. E-prints:
hep-th/9405050, q-alg/9610026, q-alg/9611003, q-alg/9611019, funct-an/9611003.
\item" [2]" Juriev D.V., An excursus into the inverse problem of representation
theory [in Russian]. Report RCMPI-95/04 [e-version: mp\_arc/96-477].
\item" [3]" Juriev D.V., Quantum projective field theory: quantum--field
analogs of Euler formulas [in Russian]. Teor.Matem.Fiz. 92(1) (1992) 172-176
[English transl.: Theor.Math.Phys. 92 (1992) 814-816; Quantum projective field
theory: quantum-field analogs of Euler-Arnold equations in projective
$G$-hypermultiplets [in Russian]. Teor.Matem.Fiz. 98(2) (1994) 220-240
[English transl.: Theor.Math.Phys. 98 (1994) 147-161]; Complex projective
geometry and quantum projective field theory [in Russian]. Teor.Matem.Fiz.
101(3) (1994) 331-348 [English transl.: Theor.Math.Phys. 101 (1994) 1387-1403];
Octonions and binocular mobilevision [in Russian]. Fundam.Prikl.Matem., to
appear; Belavkin-Kolokoltsov watch-dog effects in interactively controlled
stochactic dynamical videosystems [in Russian]. Teor.Matem.Fiz. 106(2) (1996)
333-352 [English transl.: Theor.Math.Phys. 106 (1996) 276-290]; On the
description of a class of physical interactive information systems [in
Russian]. Report RCMPI/96-05 [e-version: mp\_arc/96-459].
\item" [4]" Juriev D., The vocabulary of geometry and harmonic analysis on the
infinite-dimensional manifold $\Diff_+(S^1)/S^1$. Adv.Soviet Math. 2
(1991) 233-247.
\item" [5]" Sabinin L.V., On the nonlinear geometric algebra [in Russian]. In
``Webs and quasigroups''. Kalinin [Tver'], 1988, pp.32-37; Methods of
nonassociative algebra in differential geometry [in Russian]. Suppl.to the
Russian transl.of Kobayashi S., Nomizu K., Foundations of differential geometry.
V.1. Moscow, Nauka, 1982; Differential geometry and quasigroups [in Russian].
In ``Current problems of geometry. To the 60th Anniversary of
Acad.Yu.G.Reshetnyak''. Trans.Inst.Math. Siberian Branch Soviet Acad.Sci.,
Novosibirsk, 1984, v.14, pp.208-221; Mikheev P.O., Sabinin L.V., Quasigroups
and differential geometry. In ``Quasigroups and loops. Theory and
applications''. Berlin, Heldermann Verlag, 1990, P.357-430.
\item" [6]" Juriev D.V., Quantum conformal field theory as infinite-dimensional
noncommutative geometry [in Russian]. Uspekhi Matem.Nauk 46(4) (1991) 115-138
[English transl.: Russian Math.Surveys 46(4) (1991) 135-163].
\item" [7]" Juriev D., Infinite dimensional geometry and quantum field theory
of strings. I. Infinite-dimensional geometry of second quantized free strings.
Alg.Groups Geom. 11 (1994) 145-179 [e-version: hep-th/9403068];
II. Infinite-dimensional noncommutative geometry of self-interacting string
field. Russian J.Math.Phys. 4(3) (1996); III. Infinite-dimensional
W--differential geometry of second quantized free string. J.Geom.Phys. 16
(1995) 275-300.
\item" [8]" Taylor M., Pseudo-differential operators. B., 1974; Duistermaat J.,
Fourier integral operators. N.Y., 1973; Treves F., Introduction to
pseudo-differential and Fourier integral operators. N.Y., 1980.
\item" [9]" Shtern A.I., Quasirepresentations and pseudorepresentations [in
Russian]. Funkts.anal.i ego prilozh. 25(2) (1991) 70-73.
\item"[10]" Segal G., Unitary representations of some infinite dimensional
groups. Commun.Math. Phys. 80 (1991) 301-342; Kirillov A.A.,
Infinite-dimensional Lie groups, their orbits, in\-variants and representations.
Lect.Notes Math. 970, Springer, 1982; K\"ahler structure on K-orbits of the
group of diffeomorphisms of a circle [in Russian]. Funkts.anal.i ego prilozh.
21(2) (1987) 42-45; Kirillov A.A., Juriev D.V., K\"ahler geometry of the
infinite dimensional space $M=\Diff_+(S^1)/\Rot(S^1)$ [in Russian].
Funkts.anal.i ego prilozh. 21(4) (1987) 35-46; Kirillov A.A., Juriev D.V.,
Representations of the Virasoro algebra by the orbit method. J.Geom.Phys. 5
(1988) 351-363 [reprinted in ``Geometry and Physics. Essays in honour of
I.M.Gelfand''. Eds.S.G.Gindikin and I.M.Singer, Pitagora Editrice, Bologna
and Elsevier Sci.Publ., Amsterdam, 1991].
\item"[11]" Batalin I., Quasigroup construction and first class constraints
J.Math.Phys. 22 (1981) 1837-1850; Mikheev P.O., On the transformation loops
[in Russian]. In ``Some applications of differential geometry''. Moscow, 1985,
pp.85-93 [VINITI 4531-85Dep]; Quasigroups of transformations.
Trans.Inst.Phys.Estonian Acad.Sci. 66 (1990) 54-66.
\item"[12]" Karasev M.V., Maslov V.P., Nonlinear Poisson brackets: geometry
and quantization. Amer.Math.Soc., Providence, RI, 1991.
\item"[13]" Neretin Yu.A., Holomorphic extensions of representations of the
group of diffeomorphisms of the circle [in Russian]. Matem.Sb. 180(5) (1989)
635-657; Infinite dimensional groups, their mantles, trains and representations.
Adv.Soviet Math. 2 (1991) 103-171.
\item"[14]" Juriev D.V., On the determining of the univalence radius of a
regular function by its Taylor coefficients [in Russian].  Matem.Sb. 183(1)
(1992) 45-64; Juriev D., On the univalency of regular functions. Annali
Matem.Pura Appl.(IV) 154 (1993) 37-50.
\item"[15]" Juriev D.V., Topics in isotopic pairs and their representations
[in Russian]. Teor.Matem. Fiz. 105(1) (1995) 18-28 [English transl.:
Theor.Math.Phys. 105 (1995) 1201-1209]; Topics in isotopic pairs and their
representations. II. A general supercase [in Russian]. Teor.Ma\-tem.Fiz.
(1997), to appear.
\item"[16]" Juriev D.V., Classical and quantum dynamics of the noncanonically
coupled oscillators and Lie superalgebras. Russian J.Math.Phys. (1997), to
appear [e-version: funct-an/9409003]; On the dynamics of the noncanonically
coupled oscillators and its hidden superstructure. Preprint ESI-167
(December 1994) [e-version: solv-int/9503003].
\item"[17]" Faulkner J.R., Ferrar J.C., Simple anti-Jordan pairs.
Commun.Alg. 8(11) (1980) 993-1013.
\item"[18]" Semenov-Tian-Shanski{\v\i} M.A., What the classical r-matrix
is [in Russian]. Funkts.anal.i ego prilozh. 17(4) (1983) 17-33.
\item"[19]" Ro\v cek M., Representation theory of the nonlinear
$\operatorname{SU}(2)$ algebra. Phys. Lett.B 255 (1991) 554-557.
\item"[20]" Gerasimov A., Khoroshkin S., Lebedev D., Morozov A., Generalized
Hirota equations and representation theory. I. The case of
$\operatorname{SL}(2)$ and $\operatorname{SL}_q(2)$. E-print: hep-th/9405011
(1994); Mironov A., Quantum deformations of $\tau$--functions, bilinear
identities and representation theory. E-print: hep-th/9409190 (1994);
Kharchev S., Mironov A., Morozov A., Non-standard KP evolution and
quantum $\tau$--functions. E-print: q-alg/9501013 (1995).
\endRefs
\enddocument